\documentclass{amsart}

\usepackage{amssymb,amsfonts,latexsym}

\usepackage{pstricks}

\numberwithin{equation}{section}

\newtheorem{theorem}{Theorem}

\newtheorem{Theorem}[theorem]{Theorem}

\theoremstyle{definition}

\newtheorem{question}[theorem]{Question}

\newtheorem *{Theorem A}{Theorem A}
\newtheorem *{Theorem B}{Theorem B}
\newtheorem *{remarkA}{Remark}
\newtheorem{thm}{Theorem}[section]

\newtheorem{lem}[thm]{Lemma}

\newtheorem{cor}[thm]{Corollary}

\newtheorem{rem}[thm]{Remark}

\newcommand{\nd}{\text{and}}

\newcommand{\fr}{\text{frac}}

\newcommand{\ben}{\begin{enumerate}}

\newcommand{\een}{\end{enumerate}}

\begin{document}

\title
{A Model for Pairs of Beatty Sequences}

\author{Y. Ginosar}

\address{Department of Mathematics, University of Haifa, Haifa 31905, Israel}
\email{ginosar@math.haifa.ac.il}

\author{I. Yona}

\date{\today}


\begin{abstract}

Beatty sequences $[n \alpha_1+\beta_1]_{n\in \mathbb{Z}}$ and $[m \alpha_2+\beta_2]_{m\in \mathbb{Z}}$
are recorded by two athletes running in opposite directions in a round stadium.
This approach suggests a nice interpretation for well known partitioning criteria:
such sequences (eventually) partition the integers essentially when the athletes have the same starting point.

\end{abstract}

\maketitle

A remarkable observation due to S. Beatty
says that if $w$ is any positive irrational number, then the sequences
$\begin{array}{ccc}
1+w,& 2(1+w),& 3(1+w),...\\
1+\frac{1}{w},& 2(1+\frac{1}{w}),& 3(1+\frac{1}{w}),...
\end{array}$
contain one and only one number between each pair of consecutive positive integers.

Denoting $$\alpha_1:=1+w,\ \ \alpha_2:=1+\frac{1}{w},$$
the corresponding sequences of (floor) integer parts
$$S(\alpha_1):=[n\alpha_1]_{n\in \mathbb{N}},\ \
S(\alpha_2):=[m\alpha_2]_{m\in \mathbb{N}}$$
are called {\it Beatty sequences} with moduli $\alpha_1, \alpha_2$ respectively.

Note that the sum of $\alpha_1$ and $\alpha_2$ is equal to their product, i.e. they
satisfy
\begin{equation}\label{recip}
\frac{1}{\alpha_1}+\frac{1}{\alpha_2}=1.
\end{equation}
Beatty's result can thus be reformulated as follows:
\begin{Theorem}\label{A}(Beatty, see \cite{BSOH})
Let $\alpha_1, \alpha_2$ be two positive irrational numbers satisfying
$\frac{1}{\alpha_1}+\frac{1}{\alpha_2}=1$, then the Beatty sequences
$S(\alpha_1), S(\alpha_2)$ partition $\mathbb{N}$.
\end{Theorem}

The converse of Theorem \ref{A} is also valid:
since the density of a Beatty sequence $[n \alpha]_{n\in \mathbb{N}}$
in $\mathbb{N}$ is equal to $\frac{1}{\alpha}$, then
$S(\alpha_1), S(\alpha_2)$ partition
$\mathbb{N}$ only if $\alpha_1$ and $\alpha_2$ satisfy (\ref{recip})
(and hence there exists
a positive number $w$ such that $\alpha_1=1+w$ and $\alpha_2=1+\frac{1}{w}$).

In 1957 Th. Skolem generalized the above theorem to {\it non-homogeneous} Beatty sequences,
i.e. double infinite sequences of the form
$$S(\alpha,\beta):=[n \alpha+\beta]_{n\in \mathbb{Z}},\ \ \alpha\in \mathbb{R}^+,\ \
\beta\in\mathbb{R}.$$
When $\alpha_1$ and $\alpha_2$ are irrational,
the question is when $S(\alpha_1,\beta_1),S(\alpha_2,\beta_2)$
{\it eventually} partition $\mathbb{Z}$, that is any sufficiently large
(and any sufficiently small) integer
belongs exactly to one of the sequences.
\begin{Theorem}\label{BS}(Skolem \cite{Sk}, see also \cite{F,G})
Let $\alpha_1, \alpha_2$ be two positive irrational numbers
satisfying $\frac{1}{\alpha_1}+\frac{1}{\alpha_2}=1$, and let $\beta_1, \beta_2$
be real numbers. Then $S(\alpha_1,\beta_1),S(\alpha_2,\beta_2)$
eventually partition $\mathbb{Z}$ if and only if
\begin{equation}\label{nhc}
\frac{\beta_1}{\alpha_1}+\frac{\beta_2}{\alpha_2}\in \mathbb{Z}.
\end{equation}
Moreover, if (\ref{nhc}) holds, then $S(\alpha_1,\beta_1),S(\alpha_2,\beta_2)$
partition $\mathbb{Z}$
with an exception of, perhaps, one repeated integer $n_0$ and one missing integer $n_0-1$.
\end{Theorem}
Theorem \ref{A} clearly follows from Theorem \ref{BS}, since $0$ belongs to both sequences
$S(\alpha_1)=S(\alpha_1,0)$ and $S(\alpha_2)=S(\alpha_2,0)$
(whereas $-1$ is in neither one of them).
Hence $\{1,2,3,...\}$ ( as well as any subset of $\mathbb{Z}\setminus \{-1, 0\}$)
are disjointly covered by these sequences.

A rational analogue of Theorem \ref{BS} was established in 1969 by A.S. Fraenkel.
Obviously, if the moduli are rational,
then the notions of a partition and an eventual partition of the integers are the same.
The criterion in \cite{F} is given here in a slightly different formulation:
\begin{Theorem}\label{rathm}(Fraenkel \cite{F})
Let $r>s$ be two coprime positive integers.
Then the Beatty sequences
$S(\frac{r}{s},\beta_1)$ and $S(\frac{r}{r-s},\beta_2)$ partition $\mathbb{Z}$
if and only if
\begin{equation}\label{r-1}
[s\beta_1]+[(r-s)\beta_2]\equiv r-1\mod r.
\end{equation}
\end{Theorem}

This note suggests a perceptible approach to Beatty sequences by
interpreting $w$ as a certain ratio of speeds of two athletes running in opposite directions in a round stadium.
This interpretation yields new proofs of Theorems \ref{A}, \ref{BS} and \ref{rathm}
(in sections \ref{s1}, \ref{SkTh} and \ref{ratsec} respectively herein).
Other proofs of these theorems can be found in \cite{BB,F82,GOB,OB}.

In \S\ref{comparison}, we show that in both the rational and the irrational cases,
the partition condition says that the athletes essentially have the same starting point.

A natural question then arises:
\begin{question}\label{q}
Given two real numbers $\alpha_1$ and $\alpha_2$,
do there exist $\beta_1, \beta_2$ such that the corresponding Beatty sequences $S(\alpha_1,\beta_1)$ and $S(\alpha_2,\beta_2)$
are disjoint?
\end{question}
A complete answer to Question \ref{q} was given by R. Morikawa. It suggests an interesting notion of ``coprimeness"
of pairs of real numbers.

In \S\ref{secgamma} we review the answer to Question \ref{q} and apply again the running model to prove the following
case where $\alpha_1$ and $\alpha_2$ are irrational and of rational ratio.
\begin{Theorem}\label{gamma}(Morikawa \cite{M85})
Let $\alpha_1=r\gamma, \alpha_2=s\gamma$, where $\gamma$ is irrational and $r,s$ are coprime positive integers.
Then there exist $\beta_1, \beta_2$ such that $S(\alpha_1,\beta_1)$ and $S(\alpha_2,\beta_2)$
are disjoint if and only if $\gamma> 2$.
\end{Theorem}

An extensive bibliography on Beatty sequences and their relations to various
topics such as Sturmian words and Wythoff's game can be found in \cite{H,S,ZWS}.\\

{\bf Acknowledgement.} We thank Mota Frances, the third runner.

\section{Proof of Beatty's Theorem}\label{s1}
Assume that two athletes $X$ and $Y$ run in opposite directions in a round stadium
of length 1.
Their common starting point is denoted by $O$ and their speeds are $\frac{1}{\alpha_1}$
and $\frac{1}{\alpha_2}$ respectively.
In other words, $X$ and $Y$ complete a full round of the stadium in $\alpha_1$ and $\alpha_2$
time units respectively.
Since their speeds sum up to $1$, $X$ and $Y$ meet every time unit.
Each time one of them passes $O$,
the number of times $X$ and $Y$ have met so far is recorded.
It is easily verified that when $X$
passes $O$ for the $n$-th time, the number $[n\alpha_1]$ is recorded,
and when $Y$ passes $O$ for the $m$-th time the number $[m\alpha_2]$ is recorded.
The recorded sequences are therefore the Beatty sequences $S(\alpha_1)$ and $S(\alpha_2)$.

When $\alpha_1$ and $\alpha_2$ are irrational, then so is $w=\alpha_1-1$,
which, by (\ref{recip}), is just the ratio of these speeds. Hence
$X$ and $Y$ never meet exactly at $O$ (except at $t=0$).
Therefore, between two meetings of $X$ and $Y$, exactly one of them passes $O$.
It follows that any natural number can be uniquely expressed either as
$[n \alpha_1]_{n\in \mathbb{N}}$ or as $[m \alpha_2]_{m\in \mathbb{N}}$,
proving that $S(\alpha_1), S(\alpha_2)$ partition $\mathbb{N}$.
\qed

\section{The Non-Homogeneous Model}\label{NHM}
The above model can be fitted for non-homogeneous Beatty sequences
$S(\alpha_1,\beta_1)$ and $S(\alpha_2,\beta_2)$ with any positive real numbers
$\alpha_1,\alpha_2$
as follows:

Firstly, since for any integer $m$, the set of values of the double infinite Beatty sequences $S(\alpha,\beta)$
and $S(\alpha,\beta+m\alpha)$ are equal, then by adding or subtracting an
appropriate multiple of $\alpha_i$ to $\beta_i$, we may assume
\begin{equation}\label{beta}
0\leq \beta_i<\alpha_i
\end{equation}
for $i=1,2$.

As before, let two athletes $X$ and $Y$ run in opposite directions
in a round stadium of length 1 with speeds $\frac{1}{\alpha_1}$
and $\frac{1}{\alpha_2}$ respectively.
At time $t=0$, place $X$ and $Y$ at the points whose distances from $O$ equal
$\frac{\beta_1}{\alpha_1}$ and $\frac{\beta_2}{\alpha_2}$ respectively,
opposite their running directions.

Whenever one of the athletes passes $O$, at time $t$ say, the number $[t]$ is recorded.

It is easily verified that when $X$
passes $O$ for the $n$-th time, the number
$[(n-1)\alpha_1+\beta_1]$ is recorded,
and when $Y$ passes $O$ for the $m$-th time the number
$[(m-1)\alpha_2+\beta_2]$ is recorded.

Assume that the athletes have been recording the integers $[t]$ also while running at negative times $t<0$.
Then the running model produces the non-homogeneous Beatty sequences
$S(\alpha_1,\beta_1)$ and $S(\alpha_2,\beta_2)$.

Evidently,
these sequences are disjoint (alternatively, cover the integers) if
and only if at most (alternatively, at least) one of these athletes passes $O$ between any two consecutive integer time units.
The sequences are eventually disjoint (eventually cover the integers) if and only if the above properties
are respectively satisfied as from some $k\in \mathbb{Z}$.

\begin{rem} A set
of $n$ non-homogeneous
Beatty sequences $S(\alpha_1,\beta_1),...,S(\alpha_n,\beta_n)$ can similarly be modelled as follows:
firstly, we may once again assume (\ref{beta})
for any $1\leq i\leq n$
without changing the values of the Beatty sequences.

Next, let $X_0,X_1,...,X_{n}$ be $n+1$ athletes running in a round stadium of length 1
in the {\it same} direction
with speeds $x_0,x_1,...,x_n$ respectively. Their speeds relative to the slowest athlete $X_0$ are given by
$$x_j-x_0:=\sum_{i=1}^{j}\frac{1}{\alpha_i},\ \ 1\leq j\leq n.$$
At time $t=0$, place the $j$-th athlete at the point whose distance from $X_0$ opposite
the running direction equals
$$d_j:=\fr(\sum_{i=1}^{j}\frac{\beta_i}{\alpha_i}),\ \ 1\leq j\leq n,$$
where $\fr(u):=u-[u]$ denotes the fractional part of a real number $u$.

If at time $t$, $X_{i}$ ($1\leq i\leq n$) passes $X_{i-1}$ (which is indeed slower), then the number $[t]$ is recorded.
It is easily verified that $X_i$ records the Beatty sequence $S(\alpha_i,\beta_i)$
for any $1\leq i \leq n$.

It follows that the sequences $S(\alpha_1,\beta_1),...,S(\alpha_n,\beta_n)$ are disjoint (alternatively, cover the integers)
exactly if at most (alternatively, at least) one of the athletes $X_i$ passes $X_{i-1}$ between two consecutive integer time units.
The two athletes' model in \S \ref{NHM} can be regarded as a special instance for $n=2$,
where $X_0$ is $X$, $X_2$ is $Y$ and
where we ``accompany" the athlete $X_1$ as a steady point $O$.
\end{rem}

\section{Proof of Skolem's Theorem}\label{SkTh}
Define two domains in the stadium:
let $A$ be the set of all points whose distance
from $O$ are less than $\frac{1}{\alpha_1}$ opposite the running direction of $X$
and let $B$ be the set of all points
whose distance from $O$ are less than $\frac{1}{\alpha_2}$
opposite the running direction of $Y$.
Since (\ref{recip}) is assumed,
then the half closed domains $A$ and $B$ almost partition the stadium. Their intersection
is the point $O$, while their union misses only
the other edge point which we denote by $E$ (see figure 1).

\begin{figure}[htp]\label{fig}
\begin{pspicture}(0,-4)(0,4)
\pscircle(0,0){3}
\pscircle(0,0){2.5}
\pscustom[fillstyle=solid,fillcolor=gray]{
\psarc(0,0){3}{20}{90}
\psarcn(0,0){2.5}{90}{20}
\psline(2.8,1)}
\pscustom[fillstyle=hlines,fillcolor=black]{
\psarc(0,0){3}{90}{20}
\psarcn(0,0){2.5}{20}{90}}
\psarc{->}(0,0){2.2}{20}{90}\psarc{<-}(0,0){2.2}{20}{90}
\put(-.15,3.1){$O$}\put(1.6,2.6){$A$}\put(.9,1.5){$\frac{1}{\alpha_1}$}
\psarc{->}(0,0){2.2}{92}{18}\psarc{<-}(0,0){2.2}{92}{18}\put(-1.4,-1.5){$\frac{1}{\alpha_2}$}
\put(-2.6,-2.3){$B$}
\put(2.9,1){$E$}
\psline(-.5,-3.05)(-.45,-2.85)\put(-.77,-3.25){$Y$}
\psline(.5,-3.05)(.45,-2.85)\put(.5,-3.25){$X$}\put(-.1,-3.45){$d_0$}
\psarc{->}(0,0){3.2}{260}{280}\psarc{<-}(0,0){3.2}{260}{280}
\end{pspicture}
\caption{}
\end{figure}

At time $t=0$, $X$ and $Y$ see each other at the distance
\begin{equation}\label{d0}
d_0:=\fr(\frac{\beta_1}{\alpha_1}+\frac{\beta_2}{\alpha_2}),
\end{equation}
which, for convenience, we term as a {\it relative position $d_0$}.
Since the relative speed of the athletes is $\frac{1}{\alpha_1}+\frac{1}{\alpha_2}=1$,
they return to their relative position once every integer time units $t\in\mathbb{Z}$.
Suppose that $X$ and $Y$ are in a relative position $d_0$ at a certain time $t$
(which is an integer as shown above).
The next time they will be at the same distance is $t+1$, hence
$X$ will pass $O$ before $t+1$ precisely if it is in $A$ at the time $t$ and, similarly,
$Y$ will pass $O$ before $t+1$ precisely if it is now in $B$.

Now, assume that the Beatty sequences satisfy (\ref{nhc}). By (\ref{d0}), $X$ and $Y$
are at the same point at time $t=0$ and hence they actually meet every $t\in\mathbb{Z}$.
Each one of these meeting points is exactly in one of the domains $A$ or $B$, unless it is
on the two edges of these domains, either $O$ or $E$.
However, as before, when $\alpha_1$ and $\alpha_2$, and hence also $w$, are irrational,
$X$ and $Y$ may meet exactly at $O$ at most once.
If this happens, denote the number of this specific meeting by $n_0$.
Then $n_0$ is recorded twice, both by $X$ and by $Y$.
Clearly, in this case, the $n_0-1$-th meeting was in the other edge $E$
(which is neither in $A$ nor in $B$)
and hence the integer $n_0-1$ is not recorded.
Therefore, with this possible exception,
between two meetings of $X$ and $Y$, exactly one of them passes $O$.
It follows that any integer other than $n_0-1$ and $n_0$ can be uniquely expressed either as
$[n\alpha_1+\beta_1]$ or as $[m\alpha_2+\beta_2]$, proving that
$S(\alpha_1,\beta_1),S(\alpha_2,\beta_2)$
is a partition of $\mathbb{Z}\setminus \{n_0-1, n_0\}$.

Conversely,
suppose then that condition (\ref{nhc}) does not hold.
Then by a standard argument about the density
of irrational rotations on the unit circle (Jacobi, see \cite[Theorem 2.1]{T}),
there are infinitely many $t_i\in\mathbb{Z}$, in both time directions,
such that at time $t_i$
the athlete $X$ is in $A$ and the athlete $Y$ is in $B$.
Similarly, there are infinitely many $s_i\in\mathbb{Z}$, in both time directions,
such that at time $s_i$ $X$ is outside $A$ and $Y$ is outside $B$.
By the above interpretation of the domains $A$ and $B$, both sets
$S(\alpha_1,\beta_1)\cap S(\alpha_2,\beta_2)$ and $S(\alpha_1,\beta_1)^c\cap S(\alpha_2,\beta_2)^c$
are doubly infinite and therefore do not admit an eventual partition of $\mathbb{Z}$.
\qed

\section{Proof of Fraenkel's Theorem}\label{ratsec}

Firstly, by adding or subtracting a multiple of $\alpha_i$ to $\beta_i$ ($i=1,2$),
we keep assuming (\ref{beta}) changing neither the values
of our non-homogeneous Beatty sequences
nor condition (\ref{r-1}) in the hypothesis.

We apply once again the running model (\S\ref{NHM}) and make the following notation:
let $x_k$ and $y_k$ denote the distances of $X$ and $Y$ from $O$,
opposite the running direction of $X$, at the integer times $k\in \mathbb{Z}$.
We need the following
\begin{lem}\label{lemma}
With the above notation, the following are equivalent
\begin{enumerate}
\item $S(\frac{r}{s},\beta_1)$ and $S(\frac{r}{r-s},\beta_2)$ partition $\mathbb{Z}$.
\item for every ${k\in\mathbb{Z}}$, $x_k\in[0,\frac{s}{r})$ if and only if
$y_k\in(0,\frac{s}{r}]$.
\item for every ${k\in\mathbb{Z}}$, there exists an integer $0\leq j_k\leq r-1$ such that
$x_k\in [\frac{j_k}{r},\frac{j_k+1}{r})$ and $y_k\in (\frac{j_k}{r},\frac{j_k+1}{r}].$
\item there exists an integer $0\leq j_0\leq r-1$ such that
$x_0\in [\frac{j_0}{r},\frac{j_0+1}{r})$ and $y_0\in (\frac{j_0}{r},\frac{j_0+1}{r}].$
\end{enumerate}
\end{lem}
\noindent{\it Proof of Lemma \ref{lemma}.}
Note that the same argument as in the proof of Theorem \ref{BS}
yields a necessary and sufficient condition for the (eventual) partition of $\mathbb{Z}$
by $S(\alpha_1,\beta_1),$ and $S(\alpha_2,\beta_2)$ as follows:
for every integer time $t\in \mathbb{Z}$, one of the athletes is in $A$ if and only if
the other athlete is not in $B$.
It is easily checked that the half open segment $[0,\frac{s}{r})$ is the domain $A$,
whereas the half open segment $(0,\frac{s}{r}]$ is the complement of the domain $B$,
proving (1)$\Leftrightarrow$(2).
Next, $X$ and $Y$ return to the same relative position every integer time units.
More precisely, for every ${k\in\mathbb{Z}}$,
$$x_k-y_k\equiv d_0\mod \mathbb{Z}.$$
Now, since $s$ is prime to $r$, and since $X$ and $Y$ pass the distances $\frac{s}{r}$ and
$\frac{r-s}{r}$ respectively every integer time units,
then each of the sets $\{x_k\}_{k\in\mathbb{Z}}$ and
$\{y_k\}_{k\in\mathbb{Z}}$
divide the stadium into $r$ arcs of equal length $\frac{1}{r}$. In particular,
both sets admit exactly one member in every segment
$[\frac{j}{r},\frac{j+1}{r})$, for every integer $0\leq j\leq r-1$.
This observation implies that (2),(3) and (4) are equivalent,
and completes the proof of the lemma.
\qed

Back to the proof of Theorem \ref{rathm}, we make use of
condition (4) in Lemma \ref{lemma}.
This condition was shown to be necessary and sufficient for
partitioning the integers by $S(\frac{r}{s},\beta_1)$ and $S(\frac{r}{r-s},\beta_2)$.
Recall the initial location of $X$ and $Y$ given in \S\ref{NHM}:
$$x_0=\frac{\beta_1}{\alpha_1},\ \ y_0=1-\frac{\beta_2}{\alpha_2}.$$
Let $0\leq j_0\leq r-1$ be an integer such that $\frac{\beta_1}{\alpha_1}$ belongs to
the half open segment $[\frac{j_0}{r},\frac{j_0+1}{r})$.
Then condition (4), which is equivalent to the partitioning condition, says that $1-\frac{\beta_2}{\alpha_2}$ belongs to
the half open segment $(\frac{j_0}{r},\frac{j_0+1}{r}]$.
Equivalently,
$$j_0\leq s\beta_1<j_0+1\ \ \nd\ \ j_0< r-(r-s)\beta_2\leq j_0+1.$$
Consequently, $$[s\beta_1]+[(r-s)\beta_2]= r-1$$
is a necessary and sufficient condition for
$S(\frac{r}{s},\beta_1)$ and $S(\frac{r}{r-s},\beta_2)$ to
partition the integers under the assumption (\ref{beta}).
Relaxing this assumption, we obtain that (\ref{r-1}) is a necessary and sufficient partitioning condition.
Theorem \ref{rathm} is proven.
\qed
\section{Comparison between the cases}\label{comparison}
When the moduli $\alpha_1,\alpha_2$ are irrational and satisfy (\ref{recip}),
then Theorem \ref{BS} says that $S(\alpha_1,\beta_1)$ and $S(\alpha_2,\beta_2)$
eventually partition $\mathbb{Z}$ precisely when the athletes $X$ and $Y$ are at the very same point at $t=0$.
The rational case can be interpreted similarly:

Let $\alpha_1=\frac{r}{s}, \alpha_2=\frac{r}{r-s}\in \mathbb{Q}$.
Condition (\ref{nhc}), which says that the two athletes have a common starting point, is equivalent in this case to the condition
\begin{equation}\label{cond}
s\beta_1+(r-s)\beta_2\in r\mathbb{Z}.
\end{equation}
Next, note that the Beatty sequences
$S(\frac{r}{s},\beta_1)$ and $S(\frac{r}{s},\beta'_1)$ (with the assumption (\ref{beta}) on both)
are equal if and only if
\begin{equation}\label{betag1}
\frac{[s\beta_1]}{s}\leq \beta'_1<\frac{[s\beta_1]+1}{s}.
\end{equation}
Similarly, under the assumption (\ref{beta}),
$S(\frac{r}{r-s},\beta_2)$ and $S(\frac{r}{r-s},\beta'_2)$
are equal if and only if
\begin{equation}\label{betag2}
\frac{[(r-s)\beta_2]}{r-s}\leq \beta'_2<\frac{[(r-s)\beta_2]+1}{r-s}.
\end{equation}

Now, suppose the partitioning condition (\ref{r-1}) of Theorem \ref{rathm} holds. Let
$$\beta'_1:=\frac{[s\beta_1]+\nu}{s}$$
and
$$\beta'_2:=\frac{[(r-s)\beta_2]+1-\nu}{r-s},$$
for some $0<\nu<1$.

We have relocated the athletes such that
$$s\beta'_1+(r-s)\beta'_2\in r\mathbb{Z},$$
that is to a common starting point.

Since (\ref{betag1}) and (\ref{betag2}) are satisfied, then $S(\frac{r}{s},\beta_1)=S(\frac{r}{s},\beta'_1)$
and $S(\frac{r}{r-s},\beta_2)=S(\frac{r}{r-s},\beta'_2)$, and hence these Beatty sequences still partition the integers.

Conversely, if the common starting point condition (\ref{cond}) is satisfied,
then the partitioning condition (\ref{r-1}) holds precisely when
$s\beta_1\notin \mathbb{Z}$ (equivalently, $(r-s)\beta_2\notin\mathbb{Z}).$
This happens exactly when this common starting point is not one of the lattice points $\{\frac{j}{r}\}_{j=0}^{r-1}$.
We summarize the above discussion:
\begin{cor}
Let $r>s$ be coprime positive integers.
Suppose that the Beatty sequences $S(\frac{r}{s},\beta_1)$ and $S(\frac{r}{r-s},\beta_2)$
partition $\mathbb{Z}$. Then
$\beta_1$ and $\beta_2$ can be chosen (leaving these sequences unchanged) such that
with the interpretation of \S \ref{NHM}, the athletes $X$ and $Y$ are at the same point at $t=0$.
Conversely, suppose that the athletes are at the same point at $t=0$, with the additional demand
that this point is not one of the lattice points $\{\frac{j}{r}\}_{j=0}^{r-1}$.
Then the corresponding Beatty sequences
$S(\frac{r}{s},\beta_1)$ and $S(\frac{r}{r-s},\beta_2)$ partition $\mathbb{Z}$.
\end{cor}
\section{Disjoint Sequences}\label{secgamma}
By viewing arithmetical progressions as Beatty sequences with integral moduli,
one can formulate the well known Chinese remainder theorem as follows:
\begin{Theorem}(the Chinese remainder theorem)
Two integers $n,m$ are coprime if and only if $S(n,\beta_1)\cap S(m,\beta_2)\neq \emptyset$
for any $\beta_1, \beta_2\in\mathbb{R}$.
\end{Theorem}
Two real numbers $\alpha_1,\alpha_2$ satisfying the property that any two Beatty sequences
$S(\alpha_1,\beta_1)$ and $S(\alpha_2,\beta_2)$ intersect, may therefore be called {\it coprime} for short.

Th. Skolem gave the following necessary condition for two Beatty sets to be disjoint:
\begin{Theorem}\label{disj}(Skolem, see \cite{F})
Suppose $S(\alpha_1,\beta_1)$ and $S(\alpha_2,\beta_2)$ are disjoint. Then either
\begin{enumerate}
\item $\frac{\alpha_1}{\alpha_2}$ is rational, or
\item there exist positive integers $m,n$ such that
$$\frac{m}{\alpha_1}+\frac{n}{\alpha_2}=1,\ \ \frac{m\beta_1}{\alpha_1}+\frac{n\beta_2}{\alpha_2}\in \mathbb{Z}.$$
\end{enumerate}
\end{Theorem}
Condition (2) in Theorem \ref{disj} is also sufficient for the two sequences to be disjoint.
It takes care of integral multiples of complementary Beatty sequences (see Theorem \ref{BS}) and is therefore well understood.

The case $\alpha_1,\alpha_2\in \mathbb{Q}$ is captured by the following result known as
{\it the Japanese remainder theorem}:
\begin{Theorem}\label{JRT}(Morikawa \cite{M85}, see also \cite {SJ})
Let $\alpha_1:=\frac{p_1}{q_1}, \alpha_2:=\frac{p_2}{q_2}$ two rational numbers in reduced forms,
let $p:=(p_1,p_2), q:=(q_1,q_2),
u_1:=\frac{q_1}{q}, u_2:=\frac{q_2}{q}$. Then $\alpha_1,\alpha_2$ are coprime
if and only if there do not exist positive integers $k$ and $l$ such that
\begin{equation}\label{morik}
ku_1+lu_2=p-2u_1u_2(q-1).
\end{equation}
\end{Theorem}
Note that with the notation of the Japanese remainder theorem, when $q=1$
the existence of $k,l\in \mathbb{Z}^+$ that satisfy (\ref{morik}) is equivalent to the existence of $m,n\in \mathbb{Z}^+$
in condition (2) of Theorem \ref{disj} by the equations $$kp_1=mp,\ \ lp_2=np.$$
In particular, Theorem \ref{disj}(2) covers all the cases of $\alpha_1,\alpha_2\in\mathbb{Z}$ (i.e., where $q_1=q_2=1$).

Theorem \ref{gamma} completes the picture given in Theorems \ref{disj} and \ref{JRT}.
It deals with the remaining case, which falls under condition (1) in Theorem \ref{disj},
namely, the moduli $\alpha_1,\alpha_2$ are irrational numbers of rational ratio.

We describe the disjointness condition using the running model in \S\ref{NHM} as follows:
suppose that at time $t$ the athlete $X_1$ passes the point $O$ and hence records the integer $[t]$.
At the same time $t$, let the athlete $X_2$ be at a point whose distance from $O$
opposite the running direction is $\rho(t)$.
The disjointness condition implies that $X_2$ does not record the integer $[t]$.
This means that $X_2$ does not pass the point $O$ within the time interval $[[t],[t]+1)$.
Hence, the distance $\rho(t)$ satisfies $1-\rho(t)>\frac{\fr(t)}{\alpha_2}$
and $\rho(t)\geq \frac{1-\fr(t)}{\alpha_2}$. Consequently,
\begin{equation}\label{lattice}
\frac{1-\fr(t)}{\alpha_2}\leq \rho(t)<1-\frac{\fr(t)}{\alpha_2}.
\end{equation}
Suppose now that $\frac{\alpha_1}{\alpha_2}=\frac{r}{s}\in\mathbb{Q}$, where $r$ and $s$ are coprime positive integers
and let $\gamma:=\frac{\alpha_1}{r}=\frac{\alpha_2}{s}$.
Then $X_1$ reaches $O$ every $\alpha_1$ time units. By this time, $X_2$ passes $\frac{r}{s}$ of the stadium.
From the condition $(r,s)=1$, these steps determine a lattice $\Gamma$ of length $\frac{1}{s}$.

We are ready to prove Theorem \ref{gamma}.
This is done by applying once more the argument about density of irrational rotations on the unit circle.
If the modulus $\alpha_1$ is irrational, then for any small $\epsilon>0$,
$X_1$ passes the point $O$ at times whose fractional parts are less than $\epsilon$ as well as at times whose fractional
parts are larger than $1-\epsilon$.
Hence, (\ref{lattice}) is satisfied in this case if and only if $O$ is located exactly between two points
of $\Gamma$, and the distance $O$ to the closest lattice points is larger than the speed of $X_2$.
That is to say $\frac{1}{2s}>\frac{1}{\alpha_2}$. This completes the proof of Theorem \ref{gamma}.
$\Box$\\
\begin{rem}
Let $r,s$ be coprime positive integers.
Then with the terminology of the Chinese Remainder Theorem, Theorem \ref{gamma}
says that if $\gamma$ is irrational, then
$r\gamma$ and $s\gamma$ are coprime if and only if $\gamma< 2$.
When $\gamma\in\mathbb{Q}$,
it is not hard to see that the condition $\gamma<2$ implies that (\ref{morik}) is never satisfied.
Hence, $\gamma<2$ is sufficient for
$r\gamma$ and $s\gamma$ to be coprime .
In fact, when $\gamma$ is an integer,
the condition $\gamma<2$ is clearly necessary and sufficient for $\alpha_1=r\gamma$ and
$\alpha_2=s\gamma$ to be coprime.
\end{rem}

\end{document}